\documentclass{amsart}
\usepackage{amssymb,amsmath}
\usepackage{amsfonts}
\usepackage[all]{xy}
\usepackage{pstricks, pst-node}
\numberwithin{equation}{section}
\theoremstyle{plain}
\newtheorem{lemma}{Lemma}
\newtheorem{proposition}{Proposition}
\begin {document}
\title{About a Moduli Space of Elliptic Curves and the Golay Code $\textit{G}_{24}$}
\author{K.~M.Bugajska}
\address{Department of Mathematics and Statistics,
York University,
Toronto, ON, M3J 1P3}
\email{bugajska@yorku.ca}
\date{\today}
\begin{abstract}
We investigate algebraic structures related to triangle decompositions of a moduli space of complex tori given by the Veech curve $\textbf{T}^{*}$.  We show that these structures produce the binary error correcting Golay code $\textit{G}_{24}$.
\end{abstract}


\maketitle

\section{Intoduction}

\subsection{The Monster Group and VOA's}
 The relation between the monster group $\textbf{M}$ and the Klein modular function $\textit{J}(\tau)$ still does dot have an answer. The huge progress in this direction was made by Frenkel,Lepowski and Meurman construction ~\cite{FM88} of the vertex operaror algebra $V^{\sharp}$, i.e. the infinite dimentional graded representation of the sporadic group $\textbf{M}$. This has made the connection between $\textbf{M}$ and the modular invariant $\textit{J}$ more natural. In other words, the fact that the graded dimention of $V^{\sharp}$ is precisely equal to $\textit{J}(\tau)-744$ makes the monsterous moonshine module (the structure conjectured by Conway and Norton ~\cite{CN79} and by Thompson ~\cite{GT79}) less mysterious. It is important to notice that both objects: the monster group  $\textbf{M}$ and the vertex operator algebra $V^{\sharp}$ were developed from a sequence of structures that starts with the Golay error correcting code $\textit{G}_{24}$ and then continues with the Leech lattice $\Lambda$,  a close relative of $\textit{G}_{24}$ and determined by it.

Lepowski and Meurman also gave an alternative construction that is based on three copies of the $E_{8}$ root lattice. This construction corresponds to the fact that the code $\textit{G}_{24}$ can be obtained from copies of some (unique up to isomorphism) Hamming code that further results in the construction of $\Lambda$ as the non-orthogonal direct sum of three rescalled copies of the $E_{8}$ root lattice.

It is well known that to any even, unimodular lattice $L$ with the automorphism group $G$ we may associate a few new structures ~\cite{G06}. Thus, any such lattice $L$ allows us to produce the following objects:
\begin{itemize}
\item Thompson series $\Theta_{L}$ where  ${\Theta^{g}_{L}}={\sum{\alpha_{n}(g){q^{n}}}}$ for any element $g$ in the automorphism group $G$.
\item Thompson series $\Omega_{G}$ with  ${\Omega_{g}}={\eta_{g}(q)}$, where $\eta_{g}(q)$ is an appropriate, uniquely determined by an element $g$ product of Dedekind eta functions. The ratio $\Theta_{L}/\Omega_{G}$  forms a meromorphic Thompson series for $G$.
\item Untwisted vertex operator algebra $V_{L}$
\item A twisted vertex operator algebra $V_{L}'$.
\end{itemize} 

The real task is to construct a vertex operator algebra assoctiated to a lattice $L$ that gives the moonshine module for the group that contains $G$ or its particular extension. In our case, an even unimodular lattice is given as the Leech lattice $\Lambda$. The vertex operator algebra that forms monsterous moonshine is constructed as a $\mathbb{Z}_{2}$-orbifold obtained from the torus $\mathbb{R}^{24}/\Lambda$. More precisely, the vertex operator algebra $V^{\sharp}$ is equal to ${V^{\sharp}}={V_{1}^{+}+V_{1}^{'+}}$ where $V_{1}$ and $V'_{1}$ denote untwisted and twisted VOA's associated to the Leech lattice $\Lambda$ and the index ``+'' means that we consider only subspaces that are fixed by some involution in $\textbf{M}$.

Although both, the  investigations of Conway, Norton and Thompson (as well as many others) and the sophisticated constructions given by Frenkel, Lepowski and Meurman (and many others) are very beautiful, their connections to the Golay error correcting code $\textit{G}_{24}$, which is the origin of all structures mentioned above, is still mystierious. In other words, the question: ``What is the relation between the error correcting code $\textit{G}_{24}$ and elliptic curves '' has never been asked. We hope that in this paper we did answer this question.

\subsection{Some Teichmueller Discs and Billiards}
In our approach we will use the  following relation between some subspace of the moduli space of a concrete compact Riemann surface of genus $g>0$, say $\Sigma$, and the punctured surface $\Sigma^{*}$ itself. This surprising fact results from both: the natural hyperbolic Poincare metric on $\Sigma$ and some (singular) euclidean metric on it with appropriate symmetries. A flat structure with cone type singularities (when the genus is $g>1$) on $\Sigma$ is given by the horizontal and vertical trajectories of some holomorphic quadratic Jenkins-Strebel differential $q$ on $\Sigma$, ~\cite{M86}. Moreover, the differential $q$ determines a Teichmueller disc $\mathcal{D}_{q}(\Sigma)$ in the Teichmueller space $\mathcal{T}_{g}$. When $q$ has symmetries that produce the stabilizer group $\mathcal{G}$  of this disc in the full modular group $\mathcal{M}_{g}$ of genus $g$, then the quotient of $\mathcal{D}_{q}(\Sigma)$ by the commutator subgroup $\mathcal{G}'={[\mathcal{G},\mathcal{G}]}$ is exactly the punctures surface $\Sigma^{*}$ itself.
For example, let $\textbf{X}_{n}$ denote the hyperelliptic curve of genus $g=\frac{n-1}{2}$ that is given by the equation  $y^{2}=1-x^{n}$. The quadratic differential $q=(\frac{dx}{y})^{2}$ defines a flat metric on $\textbf{X}_{n}$ with cone type singularities when $g>1$ as well as  a Teichmueller disc $\mathcal{D}(\textbf{X}_{n})$ in$\mathcal{T}_{g}$ with the origin in $\textbf{X}_{n}$. Since $q$ determines a decomposition of the surface $\textbf{X}_n$ into $g$ vertical cylinders, the product of Dehn twists about the core loops of these cylinders forms a parabolic generator $\sigma_{n}$ of the Schwarz triangle group $\mathcal{G}_{n}=\left\langle \sigma_{n},\beta_{n}\right\rangle$, ~\cite{WV89}. The second generator $\beta_{n}$ also comes from the natural symmetries of the differential $q$.  The stabilizer of this disc in the modular group $\mathcal{M}_{g}$ is the Shwarz group $\mathcal{G}_{n}$. The compactification of he quotient of $\mathcal{D}(\textbf{X}_{n})$ by the commutator group $\mathcal{G}_{n}'$ of $\mathcal{G}_{n}$ is the Riemann surface $\textbf{X}_{n}$ again. In other words the punctured Riemann surface $\textbf{X}_{n}^{*}$ itself may be identified as a moduli space of some family of Riemann surfaces which also contains  the surface $\textbf{X}_{n}$ itself. Thus we have :
\begin{equation*}
{\textbf{X}_{n}}\in{\mathcal{D}(\textbf{X}_{n})}\hookrightarrow{\mathcal{T}_{g}}\\
\end{equation*}
 and 
\begin{equation*}
{\mathcal{D}(\textbf{X}_{n})}\rightarrow{\mathcal{D}(\textbf{X}_{n})/{\mathcal{G}_{n}'}}\cong{\textbf{X}_{n}^{*}}
\end {equation*}
  
 Pairs : [a compact Riemann surface, a holomorphic quadratic differential] that are of special importance to us are those coming from dynamical systems of billiards in rational poligons. Namely, the elastic reflections of a point-sphere particle moving with a constant velocity on the boundary of a polygon (i.e. a billiard) result in some identities. These identities glue an appropriate number of polygons together to form a compact Riemann surface that is equipped with a concrete flat metric structure whose cone type singularities occur at some vertices of glued poligons. The gluing process corresponds to the straightening out billiard trajectories. The genus of such surface is determined merely by the shape of the original polygon ~\cite{ZK76}. When all angles of the polygon have the form $\frac{\pi}{k}$ with $k\in\mathbb{Z}^{+}$, then the billiard flow is without singular points and hence the corresponding Riemann surface must be a torus. Such case is realized only for billiards in rectangles and in triangles with angles $(\frac{\pi}{3},\frac{\pi}{3},\frac{\pi}{3})$ or $(\frac{\pi}{2},\frac{\pi}{4},\frac{\pi}{4})$ or $(\frac{\pi}{2},\frac{\pi}{3},\frac{\pi}{6})$. Billiards in triangles $(\frac{\pi}{2},\frac{\pi}{n},\frac{(n-2)\pi}{2n})$ produce hyperelliptic Riemann surfaces $\textbf{X}_{n}$ together with their flat metric structure mentioned above and the group $\mathcal{G}_{n}$. The Shwarz groups $\mathcal{G}_{n}$ are conjugated in $SL_{2}\mathbb{R}$ to the Hecke triangle groups $\mathcal{H}_{n}$ generated by $S$  and by $h_{n}$ where
\begin{equation*}
S=
\begin{pmatrix}
  0&1\\
  -1&0
\end{pmatrix}
\qquad
h_{n}=
\begin{pmatrix}
  1&2cos\frac{\pi}{n}\\
  0&1
\end{pmatrix}
\end{equation*}
 
 The Hecke group $\mathcal{H}_{3}$ is exactly the modular group $\Gamma=SL_{2}\mathbb{Z}$ and the associated billiard is the  $(\frac{\pi}{2},\frac{\pi}{3},\frac{\pi}{6})$-triangle billiard which produces the torus ${\textbf{T}}\cong{\textbf{X}_{3}}$ obtained by gluing 12 such triangles together. The torus $\textbf{T}$ is the one point compactification of the quotient of the Teichmueller disc $\mathcal{D}(\textbf{T})$ by the commutator subgroup $\Gamma'=[\Gamma,\Gamma]$. In this special case the disc $\mathcal{D}(\textbf{T})$ is the whole Teichmueller space $\mathcal{T}_{1}$ and hence its quotient by the full modular group $\mathcal{M}_{1}=\Gamma$ produces the whole modular space $\mathbb{M}_{1}$ for complex tori. Since the quotient  of $\mathcal{D}(\textbf{T})$ by $\Gamma'$ produces $\textbf{T}^{*}$, each point of it represents some complex torus and points in the same orbit of the automorphisms group of $\textbf{T}^{*}$ correspond to the isomorphic tori. In a generic case, any complex torus (i.e. any point of $\mathbb{M}_{1}$) is represented by 6 points of the Veech curve   $\textbf{T}^{*}$.

 \subsection{Some Additional Symmetries}
 The deep relation between the number theory and automorphic forms for different congruence subgroups of the modular  group $\Gamma=SL_{2}\mathbb{Z}$ originates from the fact that the upper half-plane $\textit{H}$ (given by the set $\{\tau\in\mathbb{C};\mathfrak{I}\tau>0\}$)  is the moduli space of the lattices   $L_{\tau}=\mathbb{Z}+\mathbb{Z}\tau$ in $\mathbb{C}$, of complex tori $\textbf{T}_{\tau}$ (the quotient  $\mathbb{C}/L_{\tau}$) and of quadratic forms $Q_{\tau}(x,y)=\frac{1}{\mathfrak{I}\tau}|x+y\tau|^{2}$.  All these objects are interrelated to each other and any $\gamma\in\Gamma$ acting on $\textit{H}$ produces the torus $\textbf{T}_{\gamma\tau}$ that is isomorphic to $\textbf{T}_{\tau}$ and the quadratic form $Q_{\gamma\tau}(x,y)$ that is equivalent to $Q_{\tau}(x,y)$.
 
 For example, the quadratic form $Q_{\tau}(m,n)$ with $m,n\in\mathbb{Z}$ gives the lengths of closed geodesics od the torus $\textbf{T}_{\tau}$  which leads  to the geometric part of the problem (involving the counting of closed geodesics) and to the nonholomorphic Eisenstein series $E(\tau,s)$.  Its closed relative $E^{*}(\tau,s)$ given by the product $\Gamma(s)\pi^{-s}E(\tau,s)$ is related to the theta series $\Theta_{\tau}(t)$ (associated to the quadratic form $Q_{\tau}$) by the Mellin transform, etc.
 
   The similar situation occurs when $\Sigma$  is a compact Riemann surface of genus $g>1$ equipped with a (singular) flat structure determined by some    
Jenkins-Strebel differential $q$ with appropriate symmetries $\mathcal{G}$ (as mentioned in the previous subsection). In this case, if $\alpha$ is a cusp of $\Sigma^{*}\cong{\mathcal{D}_{q}(\Sigma)/\mathcal{G}'}$ then the Eisenstein-Maas series $E_{\alpha}(z,s)$ (a pure hyperbolic entity) can be expressed in terms of the `single' length spectrum of the closed leaves  in the cylinders of a flat structure corresponding to $(q,\alpha)$ (a pure flat entity). Let us notice that we are dealing not only with the hyperbolic-euclidean duality. The group $\mathcal{G}'$  (called the Veech group of $(\Sigma,q)$) not only has a finite covolume but it also has important dynamical properties associated to billiards in rational polygons. In this case we say that we have ``Veech dichotomy'' ~\cite{ZA06}. 
   
 Before we will give an outline of this paper we would like to  mention two different situations that are important for the comparison to our approach. First, we have the Millington approach (~\cite{MHM67}, ~\cite{MKAS01}) which associates to a conjugacy class of subgroups of $PSL_2\mathbb{Z}$ a transitive permutation representation of $PSL_2\mathbb{Z}$. More precisely, if $G$ is a subgroup of $\widetilde{\Gamma}=PSL_2\mathbb{Z}$ of finite index n then the decomposition of a fundamental domain of $G$ into n copies of a fundamental domains of $\widetilde{\Gamma}$ leads to the permutation group ${\Sigma(G)=\left\langle x,y\right\rangle}<S_n$. This permutation group acts transitively on the set of n letters (set of cosets) and  its generators $x; x^2=1$ and $y; y^3=1$ are the permutation images of the generators $S;\tau\rightarrow\frac{-1}{\tau}$ and  ST ($T:\tau\rightarrow\tau+1$) of the modular group $\widetilde{\Gamma}$ respectively.  If a  subgroup $G$ of index $n$ has signature $(g;e_{2},e_{3},t)$ then the permutation $x$ fixes $e_{2}$ letters , permutation $y$ fixes $e_{3}$ letters  and their product $xy$ consists of t disjoint cycles of lengths $n_{i}$ with $i=1,\ldots,t$ corresponding exactly to the cusp split $n=\sum{n_{i}}$. The different choices of fundamental domains ot the groups involved may produce transitive permutation groups whose generators $\widetilde{x},\widetilde{y}$ are necessarily simultaneously conjugate to the genarators $x,y$. For example, for $G=\Gamma(2)$ and for its quadrilateral domain ${\mathfrak{F}(\Gamma(2))}=(-1,0,1,\infty)$ we have $\mathfrak{F}(\Gamma(2))={\bigcup}_{\sigma\in\mathfrak{S}_1}{\sigma{F(\Gamma)}}$ where ${F(\Gamma)}={(i-1,\rho,i,\infty)}$ and $\mathfrak{S}_1=\{I,g,g^2,T, Tg,Tg^2\}$ with $g=ST$. When we choose the triangle ${F_{\Gamma}}=(0,\rho+1,\infty)$ as a fundamantal domain for $\widetilde{\Gamma}$ then we have $\mathfrak{F}(\Gamma(2))={\bigcup}_{\sigma\in\mathfrak{S}_2}{\sigma{F_{\Gamma}}}$  where $\mathfrak{S}_2=\{I,a,a^2,S,Sa,Sa^2\}$ and $a=TS$. It is easy to check that in both cases the generators of the transitive permutation group $\Sigma(\Gamma(2))<S_6$ are simultaneously conjugate to $x=(0 3)(1 4)(2 5)$ and to $y=(0 4 2)(3 5 1)$. Their product $xy=(01)(23)(45)$ corresponds to three cusps of width 2 each.
   On the other hand, taking the  same quadrilateral domain $(-1,0,1,\infty)$ as a fundamental domain for the index 6 subgroup $\Gamma'$ of the modular group $\widetilde{\Gamma}$  will produce  the permutation group $\Sigma(\Gamma')={\left\langle x',y'\right\rangle}$ (determined uniquely up to a simultaneous conjugation in $S_6$) with $x'=(03)(14)(25)$, $y'=(042)(153)$  and with $x'y'=(012345)$  corresponding to the single cusp of width 6.
   
  Another important situation occurs when a Riemannn surface $\textbf{X}$  is defined over an algebraic number field, or equivalently, when there exists a Belyi function $\beta:{\textbf{X}}\rightarrow{P_{1}{\mathbb{C}}}$  unbranched outside $\{0,1,\infty\}$. The triangulation $\textsl{T}_1$ of $P_{1}{\mathbb{C}}\cong{\widehat{\mathbb{C}}}$ (given by the three vertices  $0,1$ and  $\infty$;  the three edges along the line segment $\widehat{\mathbb{R}}$ joining these vertices and two triangle faces given by the upper $\textsl{H}^{+}$ and the lower $\textsl{H}^{-}$ half planes  of $\mathbb{C}\subset{\widehat{\mathbb{C}}}$ respectively) induces a triangulation $\textsl{T}={\beta^{-1}(\textsl{T}_1)}$ of $\textbf{X}$ with $\beta^{-1}(\textsl{H}^{+})$ and $\beta^{-1}(\textsl{H}^{-})$ as  open cells and with $\beta^{-1}\{0,1,\infty\}$ as vertices. (It turns out that it is more convinient to work with the Grothendieck dessins $\beta^{-1}[0,1]$ instead of with the whole $\beta^{-1}(\textsl{T}_1)$.) The monodromy group of a dessin is defined to be the monodromy group of the branching covering $\beta:{\textbf{X}}\rightarrow\widehat{\mathbb{C}}$. This is a two-generator  transitive subgroup $\left\langle g_0,g_1\right\rangle$ of the symmetric group $S_N$, where $N=deg\beta$. To find the generators $g_i, i=0,1$ we look how all edges of the bipartite graph $\beta_{-1}[0,1]$  indicent with a unique vertex in $\{\beta^{-1}(i)\}$ are permuted by rotations around their common incident vertex. The cyclic ordering around vertices $\beta^{-1}(i), i=0,1$ form the disjoint cycles of the permutations $g_i$ respectively ~\cite{GJMS97}. Two dessins are isomorphic if and only if their monodromy generators are simultaneously conjugate. Since the absolute Galois group $\textsl{Gal}(\overline{\mathbb{Q}}/\mathbb{Q})$ acts naturally on Belyi pairs $(\textbf{X},\beta)$ we have its induced action on dessins. The purpose of those investigations is to find explicitly the structure of the absolute Galois group.
  
  It is important to notice that to find the monodromy group of a dessin we may use instead of the bipartite graph $\beta^{-1}[0,1]$ the introduced earlier triangulation $\textsl{T}$ of $\textbf{X}$. Now to find the generators $g_1, g_2$ of the transitive permutation group we use the orientation of $\textbf{X}$ and rotations of positive triangles aroud appropriate vertices ~\cite{GJDS94}.
  
  We see that both, the Millington and the Belyi approaches consist of, roughly speaking, the following elements:
 \begin{itemize}
 \item Decompositions of appropriate fundamental domains into copies of a fundamantal domain of the modular group $\Gamma$. 
 
 Or 
 \item Decomposition of  Riemann surfaces over algebraic number fields into positive $(\beta^{-1}(\textsl{H}^{+}))$ and negative $(\beta^{-1}(\textsl{H}^{-}))$ triangular cells.
 
  And
 \item An action by some (external)  operations on these elements.
 \end{itemize}                                                                                                                   

  In our approach we will start  with the decompositions of fundamental domain $R=(\rho,\rho+1,\infty)$  of  $\Gamma$ into two adjacent hyperbolic $(\frac{\pi}{2},\frac{\pi}{3},0)$  triangles  $\Delta_1$ and $\Delta_2$ and we will look for algebraic structures of the operations that produce all composite triangles of the standard fundamentals domains for $\Gamma'$ out of these two   triangles.  Since we obtain that the algebraic objects related to $\mathfrak{F}'_{4 }$ are the same as those related to the quadrilaterial domain for the congruence subgroup $\Gamma(2)$ (the differences between these two domains lie in the identifications on the borders) it is necessary to introduce  also the hexagonal fundamental  domain $\mathfrak{F}_{6}'$ for $\Gamma'$. Besides, the natural projection p  of the Diag.2 maps $\mathfrak{F}'_4$ and ${\mathfrak{F}_6}'$ onto a fundamental parallelogram and the Brillouin zone of the lattice $L(\rho)$ and both of these regions are important. The reason for the consideration of both $\mathfrak{F}'_4$ and $\mathfrak{F}'_6$ is in fact much deeper and will be described in the next section.\\
   In our approach we identify the punctures torus $\textbf{T}^{*}\cong{\textit{H}/\Gamma'}$ with the Riemann surface determined by the dynamical system of the billiard in $\textbf{P}$ and we consider the group $\Gamma'$ as the Veech group of $\textbf{T}$ associated to its natural holomorphic quadratic differential.  In other words we consider $\textbf{T}^{*}$ as a Veech modular space of  complex tori (as described in the subsection $(1.2)$) instead of as (the more popular)  base curve of the elliptic modular surface $B_{{\Gamma}'}'$ (where $B_{{\Gamma}'}$ over ${\textbf{T}^{*}}\cup{\{\infty\}}$ is semistable, of arithmetic genus zero, of geometric genus one  with one singular fibre of type $I^{*}_6$, ~\cite{LL08}). Since the fibre of $B_{{\Gamma}'}'$ over the image in $\textbf{T}^{*}\cong{\textit{H}/\Gamma'}$ of a general point $\tau\in\textit{H}$ is an elliptic curve corresponding to the lattice $L_{\tau}$, the both realizations of $\textbf{T}^{*}$ (one given as the base  curve of $B_{{\Gamma}'}'$ and another given as a Veech curve) carry the same $\textsl{J}$-function. \\
    Thus, similarly as in the Millington and Belyi cases we have decompositions of fundamental domains into subdomains. However our decompositions are into triangles and hence are more subtle that the Millington ones.  The composite triangles can be viewed as coimages of the triangles $\textsl{H}^{+}$ and $\textsl{H}^{-}$ of $\widehat{\mathbb{C}}$ and hence are naturally devided into positive and negative as in Belyi approach. But in the contrary to the Millington and Grothendieck approaches our operations have totally internal nature.  They merely produce all composite triangles of $\mathfrak{F}'_4$ and $\mathfrak{F}'_6$ out of the composite modular triangles $\Delta_1$ and $\Delta_2$ of $R$.

   Now, viewing  $\textbf{T}^{*}$ as a Veech curve associated to the dynamical system of the billiard in $\textbf{P}$ determines natural bijections $\sigma$,  $\delta$ and $\kappa$ on the set of appropriate euclidean triangles. Since these triangles (Pict.3) are the images  of the hyperbolic  triangles (pict.1 and 2) respectively, we obtain the appropriate induced bijections between the hyperbolic triangle decompositions of $\mathfrak{F}'_{4}$ and $\mathfrak{F}'_{6}$.  
    
   These  algebraic operations defined on the set of our all 24 composite triangles $\Delta_{i}, i=1,\ldots,24$ (either hyperbolic ones or their, determined  by the Diag.2, euclidean images)  allow us to associate to each triangle $\Delta_{i}$  a  well defined subset $S_{i}$ of triangles. In other words each composite triangle $\Delta_{i}$ is in natural algebraic relations with some unique subset $S_i$ of triangles. When we represent these naturally arising subsets $S_{i}$ by bit strings we obtain the generating matrix for the error correcting code $\textit{G}_{24}$. Since the $\textit{J}$-function  determines the projection $J: \textbf{T}^{*}\rightarrow{Y(1)}\cong{\textit{H}/\Gamma}$ we may see this as some sort of a hidden structure associated to $J(\tau)$. (Some other  relations between algebraic structures associated to $\Gamma$ and $\Gamma'$ are given in ~\cite{KMAB10}.)
   
   As we have already mentioned in the section $1.1$, the Golay code $\textit{G}_{24}$ is the starting point for the Leech lattice and then for the monster group $\textbf{M}$. We have shown that for $\textit{G}_{24}$ to emerge we must consider a modular space given by the Veech curve $\textbf{T}^{*}$ instead of the full modular space $\mathbb{M}_{1}$. This fact has its confirmation in the ``hidden'' relations between the Leech lattice $\Lambda$ and elliptic curves with $\textit{J}$-invariant zero discovered by Harada and Lang ~\cite{HL89}.  They have shown that all five curves associated to some special conjugacy classes of the Conway group $.O$ of all automorphisms of the Leech lattice are elliptic curves that represent the Riemann surface $\textbf{T}^{*}$, that is, they all have $\textit{J}$-invariant zero. 
   
   Let us notice that if we were able to associate the Veech curve $\textbf{T}^{*}$ to a real physical object then this object would be equipped with an error correcting code (i.e. it could correct itself) and hence it would behave as some sort of a biological object.

\section{Some Algebraic Structures}
\subsection{ Fundamental Domains}
The Teichmueller space $\mathcal{T}_{1}\cong\textit{H}$ of compact  tori coincides with  $\mathcal{T}_{1,1}$ of punctured tori. Thus, any $\tau\in\textit{H}$ determines both: $\textbf{T}_{\tau}$ and $\textbf{T}_{\tau}^{*}\cong{\mathbb{C}-L_{\tau}/L_{\tau}}$. Any torus $\textbf{T}_{\tau}$  has $\mathbb{C}$ as its universal covering space (and hence it has a  flat structure), whereas $\textbf{T}^{*}_{\tau}$ is covered by the upper half-plane $\textit{H}$ with its natural Poincare hyperbolic metric. The group of deck transformations for the covering $\mathbb{C}\rightarrow\textbf{T}_{\tau}$ is isomorphic to $\mathbb{Z}^{2}$ whereas in the latter case it is an appropriate Friecke group $\mathcal{F}_{\tau}\subset{PSL_{2}(\mathbb{R})}$, ~\cite{C54}.   When $\tau=\rho={e^{\frac{{2\pi}i}{3}}}$  the Friecke group $\mathcal{F}_{\rho}$  is isomorphic to the  commutator subgroup $\Gamma'$ of the full modular group $\Gamma={SL_{2}\mathbb{Z}}$ of genus one. Its standard fundamental region $\mathfrak{F}_{4}'$ coincide with the set of elements of the standard fundamental domain $\mathfrak{F}(\Gamma(2))$ for the congruence group $\Gamma(2)$. However, the identifications on $\partial\mathfrak{F}_{4}'$ are given by the generators $A=\begin{pmatrix} 1&1\\ 1&2\end{pmatrix}$ and $B=\begin{pmatrix} 1&-1\\ -1&2\end{pmatrix}$ of $\Gamma'=\left\langle A,B\right\rangle$ and produce the punctured torus $\textbf{T}^{*}$, whereas the identifications on $\partial\mathfrak{F}(\Gamma(2))$ are given by generators $T^{2}=\begin{pmatrix} 1&2\\ 0&1\end{pmatrix}$ and $U=\begin{pmatrix} 1&0\\ 2&1\end{pmatrix}$ of $\Gamma(2)$ and produce $3$-punctured surface of genus zero.

 Let $\mathfrak{F}$ denote the quadrilateral $(-1,0,1,\infty)$ that is the underlying set for both  $\mathfrak{F}_{4}'$ and $\mathfrak{F}(\Gamma(2))$. There are a few natural decompositions of the set $\mathfrak{F}$, each of which leads to some interesting structure. First, the imaginary axis divides $\mathfrak{F}$ onto two ideal triangles (i.e. all vertices are at infinity).  The group $\Gamma$ acting on such ideal triangle produces Farrey tasselation of $\textit{H}$ which, among others, allows us to present any geodesic of $\textit{H}$ by a continuous fraction decomposition using the symbolic dynamics ~\cite{CS91}. The corresponding triangle group $\Delta^{*}_2$ (generated by the reflectios in the sides of such ideal triangle) has its assoctiated  group $\Delta(\infty,\infty,\infty)$ which is  exactly equal to $\Gamma(2)$. We will not use this fact here. Some other  decompositions which are given by $\mathfrak{F}={\mathfrak{S}_{1}F({\Gamma})}$  and $\mathfrak{F}={\mathfrak{S}_{2}F_{\Gamma}}$
 were introduced earlier.(For $\Gamma(2)$ the quadrilateral fundamental domain is the most natural one since its reflects the fact that $\Gamma(2)$ is a Fuchsian triangle group $\Delta(\infty,\infty,\infty)$ and hence its fundamental domain is given by the union of a universal triangle  and its adjacent one.)  Since the modular group $PSL_{2}\mathbb{Z}$  is also a triangle Fuchsian group $\Delta(\frac{\pi}{2},\frac{\pi}{3},0)$ both its domains $F(\Gamma)$ and $F_{\Gamma}$ are unions of two adjacent $(\frac{\pi}{2},\frac{\pi}{3},0)$-triangles and   the standard fundamental domain $R$ of $\widetilde{\Gamma}$  is the union of $\Delta_{1}=(\rho,i,\infty)$ and $\Delta_{2}=(i,\rho+1,\infty)$. However $\mathfrak{F}$ cannot be written  as  ${\bigcup_{k=1,..6}\gamma_{k}R}$. Instead  we have
 \begin{equation*}
 \mathfrak{F}={\mathfrak{S}_{1}\Delta_1}\cup{\mathfrak{S}_{2}\Delta_2}
 \end{equation*}
 Let us enumerate the composite triangles of $\mathfrak{F}$ as on Pict.1 $(\Delta_{k}\equiv{k})$.
 \begin{lemma}
 The triangle decomposition of $\mathfrak{F}$ introduced above has a natural algebraic structure of the cyclotomic coset decompositions of the quadratic residue $\mathcal{Q}$ and nonresidue $\mathcal{N}$ modulo 13 over $\mathbb{F}_{3}$.
 \end{lemma}
 \begin{proof}
  We have already introduced two sets $\mathfrak{S}_{1}$ and $\mathfrak{S}_{2}$ of coset representatives for both, for $\Gamma'$ and for $\Gamma(2)$  in $\widetilde{\Gamma}$: 
 \begin{equation}
  \mathfrak{S}_{1}={\{T,TST,T(ST)^{2}\}}\cup{\{ST,(ST)^{2},(ST)^{3}\}}
 \end{equation}
 and
 \begin{equation}
  \mathfrak{S}_{2}={\{S,STS,S(TS)^{2}\}}\cup{\{TS,(TS)^{2},(TS)^{3}\}}
 \end{equation}
 where $S=\begin{pmatrix} 0&1\\ -1&0\end{pmatrix}$, $T=\begin{pmatrix} 1&1\\ 0&1\end{pmatrix}$.   When we act on $\Delta_{1}$ by the set $\mathfrak{S}_{1}$ of transformations and on $\Delta_{2}$ by the transformations of $\mathfrak{S}_{2}$ we obtain all triangles of the decomposition. Since the set $\mathcal{Q}$ of quadratic residue in $\mathbb{F}_{13}$ is $\mathcal{Q}=\{1,4,3,12,9,10\}$ and the set $\mathcal{N}=\{2,8,6,11,5,7\}$ we see that $\mathfrak{S}_{1}(\Delta_{1})$ corresponds to $\mathcal{Q}$ and $\mathfrak{S}_{2}(\Delta_{2})$ corresponds to $\mathcal{N}$ precisely. Moreover the triangles $\Delta_{k}$ with $k\in\mathcal{Q}$ correspond to the positive triangular cells with respect to the both functions $\textsl{J}(\tau):{\mathfrak{F}'_4}\rightarrow{\mathbb{C}}$ and $\lambda(\tau):{\mathfrak{F}(\Gamma(2))}\rightarrow{\mathbb{C}}$ (i.e. their images under any of these mappings are given by the upper half-plane $\textsl{H}^{+}$ each). The images of the remaining triangles $\Delta_{k}, k\in\mathcal{N}$ are equal to the lower half-plane $\textsl{H}^{-}$ each and hence their form the negative cells. Since we have
 \begin{equation}
 \xymatrix{
 \Delta_{1}\ar[r]^T \ar@/^2pc/[rr]|3
 &\Delta_{4}\ar[r]^S \ar@/_2pc/[rr]|3
 &\Delta_{3}\ar[r]^T \ar@/^2pc/[rr]|3
 &\Delta_{12}\ar[r]^S \ar@/_2pc/[rr]|3
 &\Delta_{9}\ar[r]^T \ar@/^2pc/[rr]|3
 &\Delta_{10}\ar[r]^S \ar@/_2pc/[rr]|3
 &\Delta_{1}\ar[r]^T
 &\Delta_{4}\\
 }
 \end{equation}
 we may relate the decomposition $(2.1)$ of $\mathfrak{S}_1$ to the multiplication of $\mathcal{Q}$ by 3 and hence to the decomposition of $\mathcal{Q}$ into  disjoint union of cyclotomic cosets over $\mathbb{F}_{3}$
 \begin{equation}
 \mathcal{Q}={\{1,3,9\}}\cup{\{4,12,10\}}={\mathcal{C}_{1}}\cup{\mathcal{C}_{4}}
 \end{equation}
 
 Similarly, from the sequence
 \begin{equation}
 \xymatrix{
  \Delta_{2}\ar[r]^S \ar@/^2pc/[rr]|3
  &\Delta_{8}\ar[r]^T \ar@/_2pc/[rr]|3
  &\Delta_{6}\ar[r]^S \ar@/^2pc/[rr]|3
  &\Delta_{11}\ar[r]^T \ar@/_2pc/[rr]|3
  &\Delta_{5}\ar[r]^S \ar@/^2pc/[rr]|3
  &\Delta_{7}\ar[r]^T \ar@/_2pc/[rr]|3
  &\Delta_{2}\ar[r]^S
  &\Delta_{8}\\
  }
 \end{equation}
 we see that the decomposition $(2.2)$ of $\mathfrak{S}_{2}$ corresponds to the multiplication of $\mathcal{N}$ by 3 i.e. to the cyclotomic decomposition
 \begin{equation}
 \mathcal{N}={\{2,6,5\}}\cup{\{8,11,7\}}={\mathcal{C}_{2}}\cup{\mathcal{C}_{7}}
 \end{equation}
 \end{proof}

 From $(2.3),(2.5)$ we observe that the operation of the multiplication of a composite triangle $\Delta_k$ by 3 corresponds to the action by the transformation $a=TS$  or by the transformation $g=ST$. The transformation $a^m$ is equivalent modulo $\Gamma'$ to the transformation $g^m$ for $m=1,2,3$ (more precisely we have $ST=B(TS)=(TS)A$ and $(ST)^{2}=A^{-1}(TS)^{2}=(TS)^{2}B^{-1}$) but these transformations are not $\Gamma(2)$ equivalent ($a$ is equivalent to $g^2$ and $a^2$ is equivalent to g modulo $\Gamma(2)$). This means that our multiplication by 3 (i.e. sequences analoguous to $(2.3)$ and $(2.5)$) is well defined for $(\frac{\pi}{2},\frac{\pi}{3},0)$-triangle decomposition of any fundamental domain for  $\Gamma'$ but not for $\Gamma(2)$.
 The requirement of the quadrilateral domain $(-1,0,1,\infty)$  (for the  operation of multiplication of composite triangles by 3) for $\Gamma(2)$ is a consequence of the nonabelian quotient ${\widetilde{\Gamma}/\Gamma(2)}\cong{S_3}$.
 \vspace{10pt}
  
 \begin{pspicture}(-0.3,-2)(6.3,5)
 \rput(0,-0.2){\rnode{A}{$a_1$}}
 \rput(3,-0.2){\rnode{B}{$a_2$}}
 \rput(6,-0.2){\rnode{C}{$a_3$}}
 \rput(-0.3,5){\rnode{D}{$a_5$}}
 \rput(0.3,5){\rnode{E}{$\infty$}}
 \rput(5.7,5){\rnode{F}{$\infty$}}
 \rput(6.3,5){\rnode{G}{$a_4$}}
 \psline(0,0)(0,5)
 \psline(3,0,)(3,5)
 \psline(6,0)(6,5)
 \psarc(1.5,0){1.5}{0}{180}
 \psarc(4.5,0){1.5}{0}{180}
 \psarc(3,0){3}{0}{180}
 \psarc(0,0){3}{0}{90}
 \psarc(6,0){3}{90}{180}
 \psline(1.5,1.5)(1.5,5)
 \psline(4.5,1.5)(4.5,5)
 \rput(1.1,1.8){\rnode{a}{$\textbf{11}$}}
 \rput(0.6,2.5){\rnode{b}{9}}
 \rput(1.1,3.1){\rnode{c}{$\textbf{7}$}}
 \rput(1.9,1.8){\rnode{d}{3}}
 \rput(4.1,1.8){\rnode{e}{$\textbf{5}$}}
 \rput(4.9,1.8){\rnode{f}{12}}
 \rput(2.4,2.5){\rnode{g}{$\textbf{8}$}}
 \rput(3.6,2.5){\rnode{h}{10}}
 \rput(5.4,2.5){\rnode{i}{$\textbf{6}$}}
 \rput(1.9,3.1){\rnode{j}{1}}
 \rput(4.1,3.1){\rnode{k}{$\textbf{2}$}}
 \rput(4.9,3.1){\rnode{l}{4}}
 \rput(2.8,-1.5){\rnode{m}{Pict.1}}
 \end{pspicture}

  Contrary to this, the subgroup $\Gamma'$ is a character group of $\widetilde{\Gamma}$ and hence there is no problem for the operation of multiplication by 3 to be well defined on any of its fundamental regions. 
 
 Now let us introduce the hexagonal fundamental domain $\mathfrak{F}'_6$ for $\Gamma'$ together with its decomposition into hyperbolic $(\frac{\pi}{2},\frac{\pi}{3},0)$-triangles. The transformation $T$ acting on this set of triangles has only two orbits which coicide with the decomposition of triangles into positive and negative cells determined by the function $\textit{J}(\tau)$. So, $T$ produces all positive triangles starting from one of them and the same is true for the negative triangles. If we enumerate the positive triangles as before  by elements of $\mathcal{Q}$ and the negative triangles by elements of $\mathcal{N}$ then the transformation $T$ will correspond exactly to the multiplication by 4 in $\mathbb{F}^{*}_{13}$.  In other words, take any positive triangle as $\Delta'_1$  its adjacent to the right (negative) triangle as $\Delta'_{2}$ and the triangle numeration is given by the identification of the transformation $T$ with the  multiplication of triangles indexes by 4 modulo 13 respectively. One such enumeration of  hyperbolic $(\frac{\pi}{2},\frac{\pi}{3},0)$-triangles $\Delta_{k}'$ is shown on Pict.2.
 \begin{lemma}
 The triangle decomposition of the hexagonal fundamental region $\mathfrak{F}'_{6}$ described above carries  a natural structure of the multiplication by 4 of the quadratic residue $\mathcal{Q}$ and nonresidue $\mathcal{N}$ modulo 13.
 \end{lemma}
 \begin{proof}
 Again, the proof will be by construction. We observe that the action by $T$ on the triangles $\Delta_{i}'$ is the following:
 \begin{equation} \Delta_{1}'\stackrel{T}{\rightarrow}\Delta_{4}'\stackrel{T}{\rightarrow}\Delta_{3}'\stackrel{T}{\rightarrow}\Delta_{12}'\stackrel{T}{\rightarrow}\Delta_{9}'\stackrel{T}{\rightarrow}\Delta_{10}'\stackrel{T}{\rightarrow}\Delta_{1}'
 \end{equation}
 and
 \begin{equation}
 \Delta_{2}'\stackrel{T}{\rightarrow}\Delta_{8}'\stackrel{T}{\rightarrow}\Delta_{6}'\stackrel{T}{\rightarrow}\Delta_{11}'\stackrel{T}{\rightarrow}\Delta_{5}'\stackrel{T}{\rightarrow}\Delta_{7}'\stackrel{T}{\rightarrow}\Delta_{2}'
 \end{equation}
 Thus, after identifications of $\Delta_{i}'$ with $i$ , the operation $T$ coincides with the multiplication by $4$ of the quadratic residue $\mathcal{Q}$ i.e. to $1\stackrel{4}{\rightarrow}4\stackrel{4}{\rightarrow}3\stackrel{4}{\rightarrow}{12}\stackrel{4}{\rightarrow}9\stackrel{4}{\rightarrow}{10}\stackrel{4}{\rightarrow}1$ and to the multiplication by $4$ of the quadratic nonresidue $\mathcal{N}$ i.e.  to $2\stackrel{4}{\rightarrow}8\stackrel{4}{\rightarrow}6\stackrel{4}{\rightarrow}{11}\stackrel{4}{\rightarrow}5\stackrel{4}{\rightarrow}7\stackrel{4}{\rightarrow}2$ respectively.
 \end{proof}
 
 We notice immediately that although the sets of triangles $\{\Delta_{k},k\in{\mathcal{Q}}\}$ and $\{\Delta'_{k},k\in{\mathcal{Q}}\}$ are $\Gamma'$ equivalent we do not have their pointwise equivalence. Similarly $\{\Delta_{k},k\in{\mathcal{N}}\}$ and $\{\Delta'_{k},k\in{\mathcal{N}}\}$ are $\Gamma'$ equivalent (as the sets) but an element $\Delta_{k}$ is not necessarily $\Gamma'$ equivalent to $\Delta'_{k}$ for $k\in{\mathcal{N}}$. To understand this better suppose that the triangle $(\rho,i,\infty)$ (which is contained in both domains $\mathfrak{F}'_4$ and $\mathfrak{F}'_6$) has the same index i.e. that $\Delta_1={\Delta'_1}$. Since $T^{2}$ is $\Gamma'$-equivalent to the transformation $g^2$ we would have on the one side ${T^{2}\circ\Delta'_1}=4^2\Delta'_{1}=\Delta'_3$ whereas $T^{2}\circ{\Delta_1}$ is $\Gamma'$ equivalent to $g^{2}\circ\Delta_{1}=\Delta_9$.  So, since $\Delta'_k$ is not necessary $\Gamma'$ equivalent to $\Delta_k$ we have chosen (for $\mathfrak{F}'_6$) the first positive triangle from the right as $\Delta'_1$ and its adjacent (negative) as $\Delta'_2$. Using the transformation $T$  (identified with with the operation of the multiplication by 4 modulo 13)  we create all the remaining composite triangles of $\mathfrak{F}'_{6}$. Let $\widehat{\mathcal{T}}$ denote the set of transformations $\{T^{k},k=0,\ldots,5\}$. Thus $\mathfrak{F}'_{6}={\widehat{\mathcal{T}}\Delta'_1}\cup{\widehat{\mathcal{T}}\Delta'_2}$  with ${\widehat{\mathcal{T}}\Delta'_{1}}={\{\Delta'_{k}|k\in{\mathcal{Q}}}\}$ and ${\widehat{\mathcal{T}}\Delta'_{2}}=\{\Delta'_{k}|k\in{\mathcal{N}}\}$. 
 
 The lemmas 1 and 2  tell us that some algebraic informations contained in the sets of coset representatives $\widehat{\mathcal{T}}$ and (for example) $\mathfrak{S}_2$ of $\Gamma'$ in $PSL_{2}\mathbb{Z}$ are quite different. Thus, writing 
 \begin{equation*}
 {\widetilde{\Gamma}/\Gamma'}\cong{\left\langle{ S}\right\rangle}\times{\left\langle {g}\right\rangle}\cong{\left\langle {T}\right\rangle}mod{T}^6
 \end{equation*}
  we see that when we work with the quadrilateral domain $\mathfrak{F}'_4$ then, in fact we are dealing with the quotient $\widetilde{\Gamma}/\Gamma'$  that may be identified with the direct product of the finite subgroups of  $\widetilde{\Gamma}$ itself. However, when we consider the hexagonal fundamental domain $\mathfrak{F}'_6$ then more natural is to view the quotient as given by ${\left\langle {T}\right\rangle}mod{T^6}$, with the infinite order element ${T=Sg}\in{\widetilde{\Gamma}}$ being represented by the pair $(S,g)$.  On the other hand, the elements $a=TS$ and $g=ST$ (which are of finite orders and which are $\Gamma'$ equivalent) are responsible for the multiplication of (for example $\Delta_1$) by 3 whereas the transformation $T$ (which has infinite order in $\widetilde{\Gamma}$) is responsible for the multiplication of (for example $\Delta_1$) by 4.   The former case provides a strong relation between $\Gamma'$ and $\Gamma(2)$ whereas  the latter one  indicates strong relations between $\Gamma'$ and the groups $\Gamma_c$ and $\Gamma^{+}_{ns}(3)$, ~\cite{KMAB10}.  The genus zero subgroups $\Gamma_c$ and $\Gamma^{+}_{ns}(3)$ of $\widetilde{\Gamma}$ (which are associated to the non-split Cartan subgroups of $GL_{2}(2)$ and $GL_{2}(3)$ appropriately ~\cite{IC99}, ~\cite{BB09} ) have unique cusps of weight 2 and 3 respectively. Now, the branching schemes of the projections ${\textbf{X}'}\rightarrow{\textbf{X}_c}\rightarrow{\textbf{X}(1)}$ and ${\textbf{X}'}\rightarrow{\textbf{X}^{+}_{ns}(3)}\rightarrow{\textbf{X}(1)}$ over $\infty$ correspond exactly to the isomorphism  ${{\left\langle T\right\rangle}modT^{6}}\cong{{\left\langle T^{3}\right\rangle}\times{\left\langle T^{2}\right\rangle}modT^{6}}$. In other words we have
 \begin{equation}
 {\mathfrak{F}'_6}=(I\cup{T^3})F(\Gamma^{+}_{ns}(3))=(I\cup{T^2}\cup{T^4})F(\Gamma_c)
 \end{equation}
 
  where $F(G)$ denotes a fundamental domain of a Fuchsian group $G$ and  ${\textbf{X}'}$, ${\textbf{X}_c}$, ${\textbf{X}^{+}_{ns}(3)}$ and ${\textbf{X}(1)}$ are quotients of $\textit{H}^{*}={\textit{H}\cup\{\infty\}\cup{\mathbb{Q}}}$ by the groups $\Gamma'$, $\Gamma_c$, $\Gamma^{+}_{ns}(3)$ and $\widetilde{\Gamma}$ respectively.
 
 Moreover, if $\textit{J}_{2}$ denotes an absolute invariant of $\Gamma_c$  and  $\textit{J}_3$ is an absolute invariant for $\Gamma^{+}_{ns}(3)$ then ${Span}_{\mathbb{C}}\{\textit{J}_2,\textit{J}_3\}$ forms the underlying vector space of very important (although a nonunitary) representation  $\chi$ of $\widetilde{\Gamma}$. This representation  $\chi:{PSL_{2}\mathbb{Z}}\rightarrow{GL_{2}\mathbb{C}}$ is given by: ${\chi(S)}=I$ and by ${\chi(T)}={\begin{pmatrix} -1& 0\\ 0& \rho\end{pmatrix}}$.  Since ${\chi(S)}=I$ and since $\textit{J}_2$ and $\textit{J}_3$ live also on $\textbf{T}^{*}$ and coincide with the Weierstrass functions $\wp$ and $\wp'$  there, the vector space ${{Span}_{\mathbb{C}}\{\textit{J}_2,\textit{J}_3\}}\cong{{Span}_{\mathbb{C}}\{\wp,\wp'\}}$, (see ~\cite{KMAB10})  may be naturally viewed as the underlying space for the representation $\chi$ of the quotient $\widetilde{\Gamma}/{\Gamma}'$ only when it is realized by ${\left\langle T\right\rangle}modT^{6}$. So, it is the cusp of $\Gamma'$ that is important here and it is the hexagonal domain $\mathfrak{F}'_6$ which immediately produces the relations  $(2.9)$.
  
  Summarizing,  our approach requires the introduction of both: the quadrilateral $\mathfrak{F}'_4$ and the hexagonal $\mathfrak{F}'_6$ fundamental domains for $\Gamma'$. In both cases the subsets of the composite triangles enumerated by the elements of $\mathcal{Q}$ or of  $\mathcal{N}$ are $\Gamma'$ equivalent respectively but the particular enumeration  of the triangles of $\mathfrak{F}'_{4}$ and of $\mathfrak{F}'_{6}$ have to be independent. On the composite triangles of $\mathfrak{F}'_{4}$ we have naturally defined the operation of multiplication by 3 modulo 13 whereas the natural operations on the composite triangles of $\mathfrak{F}'_{6}$  is multiplication by 4 modulo 13. Both these operations have ``internal'' nature. To distinguish between these two ``memberships'' we have introduced the triangles $\Delta_{i}$ and $\Delta_{i}'$ respectively. From now on,  the triangles $\Delta_{i}'$ will be denoted by $\Delta_{12+i}$. An introduction of any further algebraic structure will keep this notations fixed. 
  
 In the next subsection we will use the fact that the quadrilateral domain $\mathfrak{F}'_4$ is mapped into a fundamental parallelogram of the lattice $L(\rho)$ in $\mathbb{C}$ with a vertex at the origin and with the border which forms a locus of real values of $\wp'(z,L(\rho))$. The hexagonal domain $\mathfrak{F}'_6$ of $\Gamma'$ is mapped into the Brillouin zone for the lattice $L(\rho)$ centered at $0\in{\mathbb{C}}$.  Its border forms a locus of pure imaginary values of $\wp'(z,L(\rho))$. It is well known that both of these regions are of a great importance.
 
 Now, when we consider ${\textbf{X}'}={\textbf{T}^{*}\cup{\infty}}\cong{\textit{H}^{*}}/\Gamma'$ then the image $\widehat{\Delta}$ of any  $(\frac{\pi}{2},\frac{\pi}{3},0)$-hyperbolic triangle of the tasselation of $\textit{H}$ given by $\Gamma'(\rho,i,\infty)$ must have  the double ``memberships'' which indicates whether $\textbf{T}^{*}$ is obtained from $\mathfrak{F}'_4$  or from $\mathfrak{F}'_{6}$ by appropriate  identifications on their borders.

 \subsection{Bijections Introduced by the Billiard Dynamics}
 In a general case, when $\mathcal{F}_{\tau}$ is the Friecke group corresponding to a punctured torus $\textbf{T}^{*}_{\tau}$ we have the following commutative diagram:
 \begin{equation*}
 \xymatrix{
 \textit{H}\ar[r]^{p_{\tau}} \ar[d]_{\pi_{\tau}}
 &\textit{H}{\big/}N_{\tau}\ar[d]_{p_{\tau}'} \ar@2{<->}[r]
 &\mathbb{C}-L_{\tau}\ar[d] \ar@/_2pc/[ll]^r\\
 \textbf{T}^{*}_{\tau}\ar@{<->}[r]
 &\textit{H}{\big/}{\mathcal{F}_{\tau}}\ar@2{<->}[r]
 &\mathbb{C}-L_{\tau}{\big/}L_{\tau}\\
 &Diag.1 &\\}
 \end{equation*}
 
 Here $p_{\tau}$, $p_{\tau}'$ and $\pi_{\tau}$ are the natural projections, $r$ is a conformal mapping (local inverse of $p_{\tau}$) which satisfies ${\{r,z\}}={\frac{1}{2}\wp(z,\tau)}$ with $\{\}$ denoting the Schwarzian derivative. ${N_{\tau}}=[\mathcal{F}_{\tau},\mathcal{F}_{\tau}]$ and ${\mathcal{F}_{\tau}/N_{\tau}}\cong{\mathbb{Z}^{2}}$. When $\tau$ is equal to $\rho=e^{\frac{{2\pi}i}{3}}$ then $\textbf{T}^{*}_{\rho}=\textbf{T}^{*}$ and Diag.$1$ can be extended to the following one:
 
 \begin{equation*}
 \xymatrix{
 \textit{H}\ar[r]^p \ar[rd]^\pi \ar[dd]_{\textit{J}}
 &\textit{H}{\big/}{\textit{N}} \ar[d]^{p'} \ar@2{<->}[r]
 &\mathbb{C}-L(\rho) \ar[d] \ar@/_2pc/[ll]^r\\
 & \textbf{T}^{*}\ar@2{<->}[r] \ar[dl]
 &\mathbb{C}-L(\rho){\big/}{L(\rho)}\ar[d]^{\wp'^2}\\
 \mathbb{C} & &\mathbb{C}\ar[ll]_{\textsl{T}_1}\\
 & Diag.2 &\\
 }
 \end{equation*}
 Here $\textit{N}=[\Gamma',\Gamma']$ and $\textsl{T}_{1}:z\rightarrow{z+1}$. It is known that the mapping $p:\textit{H}\rightarrow\mathbb{C}-L(\rho)$  maps the quadrilateral fundamental domain $\mathfrak{F}_{4}'$  in $\textit{H}$ onto the fundamental parallelogram $F_{4}$ in the $z$ plane $\mathbb{C}-L(\rho)$ and the hexagonal fundamental domain $\mathfrak{F}_{6}'$ onto the regular hexagon $F_{6}$ in the z-plane ~\cite{C54}. This is illustrated on the Pict.2 and on Pict.3 respectively. Vertices $a_{1}=-1, a_{2}=0, a_{3}=1$ in the extended upper half-plane $\textit{H}^{*}=\textit{H}\cup\mathbb{Q}\cup\{\infty\}$ are translates of $\infty$ by elements of $\Gamma'$ and are identified by the $\Gamma'$ generators $A$ and $B$. The mapping $p$ maps the $\Gamma'$ translates of $\infty$ into the vertices of $L(\rho)$ and the points $a_{i}'s$  of $\textit{H}^{*}$ onto the vertices of the fundamental paralellogram of the lattice $L_{0}=L(\rho)\cong{L_{\rho}}$ in $\mathbb{C}$. Moreover, $i\in\textit{H}$ and its $\Gamma$ translates are mapped onto half poins of $L(\rho)$. On $\partial{\mathfrak{F}'_4}$ the identifications are as follows: $a_{5}a_{1}\stackrel{A}{\rightarrow}a_{3}a_{2}$ with $A(a_5)=a_3$ and $A(a_1)=a_2$;   ${a_{4}a_{3}\stackrel{B}{\rightarrow}}{a_{1}a_{2}}$ with $B(a_4)=a_1$ and $B(a_3)=a_2$ (here ${a_4}={a_5}=\infty$). On $\partial{\mathfrak{F}'_{6}}$ the identifications are: ${b_2}{b_3}\stackrel{A}{\rightarrow}{b_6}{b_5}$ with $A(b_2)=b_6$, $A(b_3)=b_5$; ${b_6}{b_7}\stackrel{B}{\rightarrow}{b_4}{b_3}$ with $B(b_6)=b_4$ and $B(b_7)=b_3$;   ${b_4}{b_5}\stackrel{C}{\rightarrow}{b_8}{b_7}$ with $C(b_4)=b_8$ and $C(b_5)=b_7$ (here ${b_1}={b_9}={\infty}$ and $C=B^{-1}A^{-1}$). These all identifications induce the appropriate identifications on $\partial{F_4}$ and $\partial{F_6}$ respectively. The decompositions of $\mathfrak{F}_{4}'$ and $\mathfrak{F}_{6}'$ onto hyperbolic $(\frac{\pi}{2},\frac{\pi}{3},0)$ triangles determine the decompositions (together with their enumeration) of $F_{4}$ and $F_{6}$ into euclidean $(\frac{\pi}{2},\frac{\pi}{3},\frac{\pi}{6})$ triangles respectively. 
\\[0.5cm]
\begin{pspicture}(-3,-2)(9.3,5)
\rput(-3.2,5){\rnode{a}{$b_1$}}
\rput(-3,1.5){\rnode{b}{$b_2$}}
\rput(-1,1.5){\rnode{c}{$b_3$}}
\rput(1,1.5){\rnode{d}{$b_4$}}
\rput(3,1.5){\rnode{e}{$b_5$}}
\rput(5,1.5){\rnode{f}{$b_6$}}
\rput(7,1.5){\rnode{g}{$b_7$}}
\rput(9,1.5){\rnode{h}{$b_8$}}
\rput(9.2,5){\rnode{i}{$b_9$}}
\psline(-3,1.732)(-3,5)
\psline(-1,1.732)(-1,5)
\psline(1,1.732)(1,5)
\psline(3,1.732)(3,5)
\psline(5,1.732)(5,5)
\psline(7,1.732)(7,5)
\psline(9,1.732)(9,5)
\psline(0,0)(0,5)
\psline(2,0)(2,5)
\psline(4,0)(4,5)
\psarc(1,0){1}{0}{180}
\psarc(3,0){1}{0}{180}
\psarc(2,0){2}{60}{120}
\psarc(0,0){2}{60}{120}
\psarc(-2,0){2}{60}{120}
\psarc(4,0){2}{60}{120}
\psarc(6,0){2}{60}{120}
\psarc(8,0){2}{60}{120}
\psline(-2,2)(-2,5)
\psline(6,2)(6,5)
\psline(8,2)(8,5)
\rput(0,-0.2){\rnode{j}{$a_1$}}
\rput(2,-0.2){\rnode{k}{$a_2$}}
\rput(4,-0.2){\rnode{l}{$a_3$}}
\rput(-0.3,5){\rnode{m}{$a_5$}}
\rput(3.7,5){\rnode{n}{$a_4$}}
\rput(3,-1){\rnode{o}{Pict.2}}
\rput(-2.5,4){\rnode{A}{$\Delta'_1$}}
\rput(-1.5,4){\rnode{B}{$\Delta'_2$}}
\rput(-0.5,4){\rnode{C}{$\Delta'_4$}}
\rput(0.5,4){\rnode{D}{$\Delta'_8$}}
\rput(1.5,4){\rnode{E}{$\Delta'_3$}}
\rput(2.5,4){\rnode{F}{$\Delta'_6$}}
\rput(3.5,4){\rnode{G}{$\Delta'_{12}$}}
\rput(4.5,4){\rnode{H}{$\Delta'_{11}$}}
\rput(5.5,4){\rnode{I}{$\Delta'_9$}}
\rput(6.5,4){\rnode{J}{$\Delta'_5$}}
\rput(7.5,4){\rnode{K}{$\Delta'_{10}$}}
\rput(8.5,4){\rnode{L}{$\Delta'_7$}}
\end{pspicture}
 
 Now, the embeddings of $F_{4}$ and $F_{6}$ into $\mathbb{C}$ and their decompositions into $(\frac{\pi}{2},\frac{\pi}{3},\frac{\pi}{6})$ euclidean triangles correspond to the tasselation of $\mathbb{R}^{2}$ by the copies of the triangle billiard in the triangle $\textbf{P}$. The process of gluing together $12$ copies of the triangle $\textbf{P}$ yields two (equivalent) tasselations of $\mathbb{R}^{2}$ given by the period parallelogram $F_{4}=(a_{1},a_{2},a_{3},a_{4})$  or by the regular hexagon $F_{6}=(b_{2},b_{3},\ldots,b_{7})$ respectively (see Pict.3).
 
 \begin{pspicture}(-2,-4.5)(4,3.6)
 \pspolygon(-1.732,-1)(-1.732,1)(0,2)(1.732,1)(1.732,-1)(0,-2)
 \pspolygon(0,0)(1.732,3)(3.464,0)(1.732,-3)
 \rput(-2,-1){\rnode{A}{$b_2$}}
 \rput(-2,1){\rnode{B}{$b_7$}}
 \rput(0,2.2){\rnode{C}{$b_6$}}
 \rput(2,1){\rnode{D}{$b_5$}}
 \rput(2,-1){\rnode{E}{$b_4$}}
 \rput(0,-2.2){\rnode{F}{$b_3$}}
 \rput(-0.3,0){\rnode{G}{P}}
 \rput(1.732,3.2){\rnode{H}{$a_3$}}
 \rput(3.8,0){\rnode{I}{$a_2$}}
 \rput(1.732,-3.2){\rnode{J}{$a_1$}}
 \rput(0.2,-4.2){\rnode{K}{Pict.3}}
 \rput(2.2,-4.2){\rnode{L}{$(P=a_4=b_9)$}}
 \end{pspicture}

 The process of straighting out billiard's trajectories corresponds to reflections of the triangle $\textbf{P}$ (or its copies) about its edges. In other words, the projection of any line $l\subset\mathbb{R}^{2}$ onto $\textbf{P}$ forms a billiard trajectory, whereas the projection of $l$ onto $F_{4}$ or $F_{6}$ (using appropriate identifications on the borders) is a flat geodesics in $\textbf{X}'$. When we start with $F_{4}$ then the dynamics of the billiard in $\textbf{P}$   corresponds to reflections of the composite triangles about the border of $F_{4}$ and hence to the reflections of the euclidean $(\frac{\pi}{2},\frac{\pi}{3},\frac{\pi}{6})$ triangles in their longer perpendicular. If we start with $F_{6}$ then the ``billiard identifications'' at its border correspond to te reflections of a composite triangle in its shorter perpendicular. Hence, the billiard dynamics introduces the following bijections on the sets of the composite triangles of $F_{4}$ and $F_{6}$ respectively:
 \begin{itemize}
 \item The bijection $\sigma$, which corresponds to the reflections of $(\frac{\pi}{2},\frac{\pi}{3},\frac{\pi}{6})$ triangles in their longer perpendicular.
 \item The bijection $\delta$ that discribes the reflections of the composite triangles in their shorter perpendicular.
 \item The bijection $\kappa$, the natural bijection between the sets of composite triangles of $F_{4}$ and $F_{6}$. 
 \end{itemize}
More precisely, the action of $\sigma$ is given by (we use the notation $k$ for $\Delta_{k}$):
\begin{align*}
&1\rightarrow8\rightarrow1      &13\rightarrow23\rightarrow13\\
&2\rightarrow10\rightarrow2     &14\rightarrow24\rightarrow14\\
&3\rightarrow11\rightarrow3     &16\rightarrow17\rightarrow16\\
&4\rightarrow12\rightarrow4     &20\rightarrow21\rightarrow20\\
&5\rightarrow12\rightarrow5     &15\rightarrow19\rightarrow15\\
&7\rightarrow9\rightarrow7      &18\rightarrow22\rightarrow18\\
\end{align*}

Similarly, for the mapping $\delta$ we have:
\begin{align*}
 &13\rightarrow14\rightarrow13    &4\rightarrow11\rightarrow4\\
 &15\rightarrow18\rightarrow15    &5\rightarrow9\rightarrow5\\
 &16\rightarrow20\rightarrow16    &6\rightarrow3\rightarrow6\\
 &17\rightarrow21\rightarrow17    &7\rightarrow12\rightarrow7\\
 &19\rightarrow22\rightarrow19    &8\rightarrow10\rightarrow8\\
 &23\rightarrow24\rightarrow23    &1\rightarrow2\rightarrow1\\
\end{align*}

 Of course we have the analogous bijections (also denoted by $\sigma$ and $\delta$) on the set of the composite triangles of $\mathfrak{F}'_{4}$ and $\mathfrak{F}'_{6}$. In other words, the above relations may be read from the reflections of the composite hyperbolic triangles in the  borders $\partial\mathfrak{F}'_{4}$ and $\partial\mathfrak{F}'_{6}$ using the borders identifications given by $A$ and $B$  or by $A$, $B$ and $C=B^{-1}A^{-1}$ respectively. However, reflections are natural for billiards and do not belong to the hyperbolic ``picture'' of $\textbf{T}^{*}$ (they belong to the Veech ``picture''instead).

The last bijection $\kappa$ [which also may be read from either the hyperbolic, Pict.1 and Pict.2, or from the billiard relations] is the following:
\begin{align*}
   &1\leftrightarrow15     &7\leftrightarrow20\\
   &2\leftrightarrow18     &8\leftrightarrow19\\
   &3\leftrightarrow13     &9\leftrightarrow21\\
   &4\leftrightarrow24     &10\leftrightarrow22\\
   &5\leftrightarrow17     &11\leftrightarrow23\\
   &6\leftrightarrow14     &12\leftrightarrow16\\
\end{align*}
Let $\Omega$ denote the set of $24$ triangles $\Delta_{i}$'s such that $\bigcup_{i=1}^{12}\Delta_{i}=F_{4}\cong{\mathfrak{F}'_{4}}$ and $\bigcup_{i=13}^{24}\Delta_{i}=F_{6}\cong{\mathfrak{F}'_{6}}$.

\section{Moduli Space and the Golay Code}

\subsection{Structures on the Set $\Omega$}
We have seen that the hyperbolic structures on $\mathfrak{F}'_{4}$  and on $\mathfrak{F}'_{6}$ induces the following:
\begin{itemize}
\item On the set of composite triangles of $F_{4}$ the natural operation of multiplication by $3$ (equivalently, the cyclotomic coset structure of the quadratic residue and nonresidue in $\mathbb{F}_{13}$ over $\mathbb{F}_{3}$).
\item On the set of composite triangles of the regular hexagon $F_{6}$ the natural operation of multiplication by $4$.
\item The natural bijection $\kappa$ between the sets of the composite triangles of these two regions. 
\end{itemize}
On the other hand the flat structure on $\textbf{T}^{*}$ determined by the dynamical system of the billiard produces the mappings $\sigma$ and $\delta$ defined in the previous subsection. When a $(\frac{\pi}{2},\frac{\pi}{3},\frac{\pi}{6})$-triangle has the $F_4$ membership (i.e. for $\Delta_i$ with $i=1,\ldots,12$) the natural operation is the multiplication by 3 and the natural mapping is $\sigma$ (since it  is introduced by the billiard reflections in the border $\partial{F}_4$).  When a triangle has the $F_6$ membership (i.e. for $\Delta_i$ with $i=13,\ldots,24$) the natural operation is the multiplication by 4 and the natural  mapping is $\delta$ (corresponding to the billiard reflections in $\partial{F}_6$). Thus, on the set $\Omega=\{\Delta_{i}| i=1,\ldots,24$ we have rather rich algebraic properties. They allow us to construct natural correspondences which to each ${\Delta_{i}}\equiv{i}$ associate a subset $S_{i}\subset\Omega$ in the following way: For $i=1,\ldots,12$
\begin{equation}
\xymatrix{
[3i]&\textbf{i}\ar[l]_{3} \ar[r]^{\sigma}\ar@{<->}[d]^{\kappa}&\sigma(i)\ar[r]^{3} \ar@{<->}[d]^{\kappa}
&[3\sigma(i)]\\
4\circ{\kappa(i)}&{\kappa(i)}\ar[l]_{4} \ar[r]^{\sigma} &{\sigma{\kappa(i)}}\ar[r]^{4}
&4\circ{\sigma{\kappa(i)}}\\} 
\end{equation}
Here the squere bracket denotes the congruence modulo 13, 3 and 4 denote the multiplications by 3  and by 4  modulo 13 respectively. Moreover, since we have to make appropriate adjustments caused by the  identification of $\Delta'_{i}$ with $\Delta_{i+12}$, we  introduce     
\begin{equation*}
     4\circ{\kappa(i)}:=[4(\kappa(i)-12)]+12
\end{equation*}
and 
\begin{equation*}
4\circ{\sigma{\kappa(i)}}:=[4(\sigma{\kappa(i)}-12)]+12
\end{equation*}

 Notice that instead of using the language of quadratic residue and nonresidue of $\mathbb{F}_{13}$, cyclotomic cosets and so on we could consider appropriate sequenses of transformations $S$ and $T$ acting on the hyperbolic triangle decompositions of $\mathfrak{F}'_{4}$ and $\mathfrak{F}'_{6}$ appropriately. However, these elementary number theory representations make our formulae much more transparent and elegant. For example, for $i=1$ we have
 \begin{equation}
 \xymatrix{
 3&\textbf{1}\ar[l]_{3} \ar[r]^{\sigma} \ar@{<->}[d]^{\kappa} &8\ar[r]^{3} \ar@{<->}[d]^{\kappa}&11\\
 24&15\ar[l]_{4} \ar[r]^{\sigma} &19\ar[r]^{4} &14\\}
 \end{equation}. 
 
  Hence, for the triangle $\Delta_1$ the associeted subset $S_1$ of the set $\Omega$ is given by $\{1,8,3,11,15,19,24,14\}$. Sorting the elements of each set $S_k$   in their increasing order we obtain the following:
 \begin{align*}
    &S_{1}=\{1,3,8,11,14,15,19,24\}=S_{8}\\
    &S_{2}=\{2,4,6,10,13,18,22,23\}=S_{10}\\
    &S_{3}=\{3,7,9,11,13,16,17,23\}=S_{11}\\
    &S_{4}=\{4,5,6,12,14,20,21,24\}=S_{6}\\
    &S_{5}=\{2,5,10,12,15,16,17,19\}=S_{12}\\
    &S_{7}=\{1,7,8,9,18,20,21,22\}=S_{9}\\
 \end{align*}
 Since the bijections $\sigma$, $\delta$ and $\kappa$ satisfy:
 \begin{equation}
   \kappa{\sigma}=\sigma{\kappa} \qquad and \qquad  \kappa{\delta}=\delta{\kappa}
 \end{equation}
 we see immediately that $S_{i}=S_{\sigma(i)}$.  Similarly, to each triangle $\Delta_{j}$ with $j=13,\ldots,24$ we associate a subset $S_{j}$ according to:
 
 \begin{equation}
 \xymatrix{
 4\circ{j} &\textbf{j}\ar[l]_{4} \ar[r]^{\delta} \ar@{<->}[d]^{\kappa}
 &\delta(j)\ar[r]^{4} \ar@{<->}[d]^{\kappa} &4\circ{\delta(j)}\\
 [3\kappa(j)] &\kappa(j)\ar[l]_{3} \ar[r]^{\delta} &\kappa{\delta(j)} \ar[r]^{3} &[3\kappa{\delta(j)}]\\}
 \end{equation}
 where as above $4\circ{j}:=[4(j-12)]+12$ and $4\circ{\delta(j)}:=[4(\delta(j)-12)]+12$ correct our notation. Thus, for example for $j=13$ we have:
 
 \begin{equation}
 \xymatrix{
 16 &\textbf{13}\ar[l]_{4} \ar@{<->}[d]^{\kappa} \ar[r]^{\delta} &14\ar@{<->}[d]^{\kappa} \ar[r]^{4} &20\\
 9 &3\ar[l]_{3} \ar[r]^{\delta} &6\ar[r]^{3} &5\\}
 \end{equation}
 Again we have $S_{j}=S_{\delta(j)}$ for $j=13,\ldots,24$ and sorting their elements by the increasing order the explicit forms of these subsets of $\Omega$ are:
 \begin{align*}
 &S_{13}=\{3,5,6,9,13,14,16,20\}=S_{14}\\
 &S_{15}=\{1,2,3,6,15,18,23,24\}=S_{18}\\
 &S_{16}=\{7,8,10,12,15,16,18,20\}=S_{20}\\
 &S_{17}=\{1,2,5,9,17,19,21,22\}=S_{21}\\
 &S_{19}=\{4,8,10,11,13,14,19,22\}=S_{22}\\
 &S_{23}=\{4,7,11,12,17,21,23,24\}=S_{24}\\
 \end{align*}
 Now, each subset $S_{i}$, $1=1,\ldots,24$ may be represented by a bit string of the lenght $24$. Passing to the bit string description we immediately notice that, for example  $S_{7}$ is equal to the sum of the strings $S_{1}+S_{2}+S_{3}+S_{4}+S_{5}$  and that $S_{23}$ is equal to the sum $S_{13}+S_{15}+S_{16}+S_{17}+S_{19}$. This is not surprising. From the correspondences $\Delta_{i}\rightarrow{S_{i}}$ it is obvious that the association $\Delta_{7}\rightarrow{S_{7}}$ as well as $\Delta_{23}\rightarrow{S_{23}}$  are already described by the remaining correspondences. However, we still did not involve the explicit correspondences between  the cyclotomic coset decompositions and the composite triangles.
 
 So,to have the full description of algebraic structures associated to the Teichmueller disc $\mathcal{D}(\textbf{T})\cong{\textit{H}}$ with symmetries determined by the appropriate holomorphic quadratic  differential (equivalently by the dynamical system of the billiard in $\textbf{P}$) we must also consider the following subsets of $\Omega$:
 \begin{equation}
 R_{6}:=\mathcal{C}_{1}\cup{\mathcal{C}'_{4}}\cup{\mathcal{C}'_{2}}\cup{\mathcal{C}_{7}}
 \end{equation}
 and
 \begin{equation}
 R_{12}:=\mathcal{C}'_{1}\cup{\mathcal{C}_{4}}\cup{\mathcal{C}_{2}}\cup{\mathcal{C}'_{7}}
 \end{equation}
 
 Here $\mathcal{C}_{1}=\{1,3,9\}$ and $\mathcal{C}_{4}=\{4,12,10\}$ are cyclotomic cosets of $\mathcal{Q}$ over $\mathbb{F}_{3}$ whereas  $\mathcal{C}_{2}=\{2,6,5\}$ and $\mathcal{C}_{7}=\{8,11,7\}$ are the cyclotomic cosets of $\mathcal{N}$.  The $\kappa$-image of a coset $\mathcal{C}_{l}$ is denoted by $\mathcal{C}'_{l}$.
 
 Geometrically, the subsets $R_{6}$ and $R_{12}$ are associated to the decomposition of $\mathfrak{F}'_{4}$ onto two ideal triangles. More precisely, $R_{6}$ represents all composite triangles of the ideal hyperbolic triangle with the vertices at $(\infty,-1,0)$ and the $\kappa$-images of the composite triangles of the other half of $\mathfrak{F}'_{4}$ (i.e. of the ideal triangle $(0,1,\infty)\subset{\mathfrak{F}'_{4}}$). The subset $R_{12}$ represents vice versa, all composite triangles of $(0,1,\infty)$ and the $\kappa$-images of $(\infty,-1,0)$. We may say, that similarly as the necessity of the involving of the hexagonal region for $\Gamma'$  is connected to the relation between the groups $\Gamma'$ and the mentioned earlier $\Gamma^{+}_{ns}(3)$ and $\Gamma_c$ the subsets given by the formulas $3.6$ and $3.7$ are connected to the relations between $\Gamma'$ and $\Gamma(2)$.
 
Our original enumeration of the composite hyperbolic $(\frac{\pi}{2},\frac{\pi}{3},0)$ triangles could be quite different. Although the decompositions into positive ( corespondending to the set of elements of $\mathcal{Q}$) triangles and negative (corresponding to the set of elements of $\mathcal{N}$) triangles are fixed (by the values of  the modular invariant $\textit{J}(\tau)$) we could take any positive triangle in $\mathfrak{F}'_4$ as $\Delta_1$ and its adjacent as $\Delta_2$ and create the remaining composite triangles of the quadrilateral using $(2.3)$ and $(2.5)$. Similarly any positive triangle of $\mathfrak{F}'_6$ and its adjacent (which is of course negative) can be taken as ${\Delta'_1}={\Delta_{13}}$ and ${\Delta'_2}={\Delta_{14}}$ respectively and the remaining composite triangles can be obtain from them using $(2.7)$ and $(2.8)$.  This means that our relations $\Delta_{i}\rightarrow{S_{i}}$ could look quite differently. However, the algebraic strucrures which are determined by the sequences $(2.3)$, $(2.5)$, $(2.7)$ and $(2.8)$ as well as by the properties of  the mappings $\sigma$, $\delta$ and $\kappa$  would be exactly the same,   merely described by a different notation.
 
 Summarizing, all natural algebraic properties that may be associated to $\textbf{T}^{*}$ seen as a  Veech moduli space of compact complex tori and related to its hyperbolic and euclidean structures, are totally described by the following  set of subsets of $\Omega$:
 \begin{equation}
 {\{S_{1},S_{2},S_{3},S_{4},S_{5},S_{13},S_{15},S_{17},S_{19}\}}\cup{\{R_{6},R_{12}\}}
 \end{equation}
  
 \subsection{Moduli Space and the Golay Code}
 In the previous subsections we have obtained same  algebraic structures that uniquely describe the nature of the  triangle decompositions of the quadrilateral and the hexagonal fundamental regions of $\Gamma'$ together with their relations to the billiard's dynamics in $\textbf{P}$. However,  instead of working with these two fundamental domains we could choose only one of them and consider each of its composite triangle twice (i.e. give each of them two, independent indices) depending whether it undergoes the operation of multiplication by $3$ and the operation $\sigma$ or it undergoes the operation of multiplication by $4$ and by $\delta$. Thus each composite triangle of any fundamental domain for $\Gamma'$ would have two (not necessarily equal) labels connected by the obvious bijection $\kappa$. Of course, the algebraic structure which describes all algebraic-geometric relations between these $24$ (single) labelled triangles is again given by the set of the form of $(3.8)$.
 
 \begin{proposition}
 The algebraic structure mentioned above and describing all algebraic informations given by the triangle decompositions of the fundamental parallelogram  and the regular hexagon  (viewed as the   $p$-images of the domains $\mathfrak{F}'_{4}$ and $\mathfrak{F}'_{6}$ of $\Gamma'$) is exactly the one given by the error correcting Golay code $\textit{G}_{24}$.
 \end{proposition}
 \begin{proof}
 First, let us rename our subsets $S_{i}$'s as follows: $R_{i}:=S_{i}$ for $i=1,\ldots,5$. Let $R_{7}:=S_{13}$, $R_{8}:=S_{15}$, $R_{9}:=S_{16}$, $R_{10}:=S_{17}$ and let $R_{11}:=S_{19}$. Now it is the set $\{R_{i}|i=1,\ldots,12\}$ that contains all informations about the algebraic properties carried out by the triangle decompositions of the moduli space $\textbf{T}^{*}$. Since each subset $R_{i}$ corresponds to a bit string we may construct a matrix $\textbf{G}$ whose i-th row is given by the bit string  of lenght $24$ determined by $R_{i}$.
 \begin{equation}
 \left(
 \begin{array}{cccccccccccccccccccccccc}
 1&0&1&0&0&0&0&1&0&0&1&0&0&1&1&0&0&0&1&0&0&0&0&1\\
 0&1&0&1&0&1&0&0&0&1&0&0&1&0&0&0&0&1&0&0&0&1&1&0\\
 0&0&1&0&0&0&1&0&1&0&1&0&1&0&0&1&1&0&0&0&0&0&1&0\\
 0&0&0&1&1&1&0&0&0&0&0&1&0&1&0&0&0&0&0&1&1&0&0&1\\
 0&1&0&0&1&0&0&0&0&1&0&1&0&0&1&1&1&0&1&0&0&0&0&0\\
 1&0&1&0&0&0&1&1&1&0&1&0&0&1&0&1&1&1&0&0&0&1&0&1\\
 0&0&1&0&1&1&0&0&1&0&0&0&1&1&0&1&0&0&0&1&0&0&0&0\\
 1&1&1&0&0&1&0&0&0&0&0&0&0&0&1&0&0&1&0&0&0&0&1&1\\
 0&0&0&0&0&0&1&1&0&1&0&1&0&0&1&1&0&1&0&1&0&0&0&0\\
 1&1&0&0&1&0&0&0&1&0&0&0&0&0&0&0&1&0&1&0&1&1&0&0\\
 0&0&0&1&0&0&0&1&0&1&1&0&1&1&0&0&0&0&1&0&0&1&0&0\\
 0&1&0&1&1&1&0&0&0&1&0&1&1&0&1&0&0&0&1&1&1&0&1&0\\
 \end{array}
 \right)
 \end{equation}
 These bit strings span $12$-dimentional subspace of the space $\mathbb{F}_{2}^{24}$. It was checked (on a computer) that the weight distribution is given exactly by
 \begin{equation}
                     1+759q^{8}+2576q^{12}+759q^{16}+q^{24} 
 \end{equation}
 We recall that the weight of a bit string is given by the number of $1$'s . Now, using well known theorem (see ~\cite{MS77}) we obtain that our set $\{R_{1},\ldots,R_{12}\}$  determines unique, up to isomorphism, error correcting binary Golay code $\textit{G}_{24}$.  
 \end{proof}
 
 We notice that our generating matrix $\textbf{G}$ for the binary code $\textit{G}_{24}$ has quite different form that the standard one. The reason for this is that usually to construct some code we are trying to find a nice matrix (i.e. a subspace of $\mathbb{F}_{k}^{n}$) whereas here our matrix comes from the algebraic structures naturally arising from the (doubly  indexed) triangle decompositions of the moduli space (given by the Veech space $\textbf{T}^{*}$) of compact complex tori. Since the modular invariant $\textit{J}(\tau)$ lives  on this moduli space and  the projection $\textit{J}:{\textbf{T}^{*}}\rightarrow{Y(1)}$ determines the decompositions of $\textbf{T}^{*}$ into all positive ($\mathcal{Q}$) and negative ($\mathcal{N}$) triangle cells we have obtained a sort of a hidden structure  which is given by $\textit{G}_{24}$ and which is associated to $\textit{J}(\tau)$. By an appropriate change of the enumeration of our triangles we may get the standard form of the generating matrix for the Golay code $\textit{G}_{24}$.

\end{document}